# OVERSHOOTS AND UNDERSHOOTS OF LÉVY PROCESSES

By R. A. Doney and A. E. Kyprianou

*University of Manchester and Heriot Watt University*

We obtain a new fluctuation identity for a general Lévy process giving a quintuple law describing the time of first passage, the time of the last maximum before first passage, the overshoot, the undershoot and the undershoot of the last maximum. With the help of this identity, we revisit the results of Klüppelberg, Kyprianou and Maller [*Ann. Appl. Probab.* **14** (2004) 1766–1801] concerning asymptotic overshoot distribution of a particular class of Lévy processes with semi-heavy tails and refine some of their main conclusions. In particular, we explain how different types of first passage contribute to the form of the asymptotic overshoot distribution established in the aforementioned paper. Applications in insurance mathematics are noted with emphasis on the case that the underlying Lévy process is spectrally one sided.

**1. Lévy processes and ladder processes.** This paper concerns overshoots and undershoots of Lévy processes at first upwards passage of a constant boundary. We will therefore begin by introducing some necessary but standard notation.

In the sequel $X$ will always denote a Lévy process defined on the filtered space $(\Omega, \mathcal{F}, \mathbb{F}, P)$, where the filtration $\mathbb{F} = \{\mathcal{F}_t : t \geq 0\}$ is assumed to satisfy the usual assumptions of right continuity and completion. Its characteristic exponent will be given by $\Psi(\theta) := -\log E(e^{i\theta X_1})$ and its jump measure by $\Pi_X$. We will work with the probabilities $\{P_x : x \in \mathbb{R}\}$ such that $P_x(X_0 = x) = 1$ and $P_0 = P$. The probabilities $\{\widehat{P}_x : x \in \mathbb{R}\}$ will be defined in a similar sense for the dual process, $-X$.

Denote by $\{(L_t^{-1}, H_t) : t \geq 0\}$ and $\{(\widehat{L}_t^{-1}, \widehat{H}_t) : t \geq 0\}$ the (possibly killed) bivariate subordinators representing the ascending and descending ladder processes. Denote by $\kappa(\alpha, \beta)$ and $\widehat{\kappa}(\alpha, \beta)$ their joint Laplace exponents for

$\alpha, \beta \geq 0$. For convenience, we will write

$$\kappa(0, \beta) = q + \xi(\beta) = q + c\beta + \int_{(0,\infty)} (1 - e^{-\beta x}) \Pi_H(dx),$$

where $q \geq 0$ is the killing rate of $H$ so that $q > 0$ if and only if $\lim_{t \uparrow \infty} X_t = -\infty$, $c \geq 0$ is the drift of $H$ and $\Pi_H$ is its jump measure. The quantity $\xi$ is a true subordinator Laplace exponent. Similar notation will also be used for $\widehat{\kappa}(0, \beta)$ by replacing $q$, $\xi$, $c$ and $\Pi_H$ by $\widehat{q}$, $\widehat{\xi}$, $\widehat{c}$ and $\Pi_{\widehat{H}}$. Note that, when $q > 0$, we necessarily have $\widehat{q} = 0$.

Associated with the ascending and descending ladder processes are the bivariate renewal functions $\mathcal{U}$ and $\widehat{\mathcal{U}}$. The former is defined by

$$\mathcal{U}(dx, ds) = \int_0^\infty dt \cdot P(H_t \in dx, L_t^{-1} \in ds)$$

and taking double Laplace transforms shows that

$$(1) \qquad \int_0^\infty \int_0^\infty e^{-\beta x - \alpha s} \mathcal{U}(dx, ds) = \frac{1}{\kappa(\alpha, \beta)} \qquad \text{for } \alpha, \beta \geq 0,$$

with a similar definition and relation holding for $\widehat{\mathcal{U}}$. These bivariate renewal measures are essentially the Green's measures of the ascending and descending ladder processes. By $U(dx)$ and $\widehat{U}(dx)$, we will denote the marginal measures $\mathcal{U}(dx, [0, \infty))$ and $\widehat{\mathcal{U}}(dx, [0, \infty))$, respectively. Note that local time at the maximum is defined only up to a multiplicative constant. For this reason, the exponent $\kappa$ can only be defined up to a multiplicative constant and, hence, the same is true of the measure $\mathcal{U}$ (and then obviously this argument applies to $\widehat{\mathcal{U}}$).

Let

$$\overline{X}_t := \sup_{u \leq t} X_u \quad \text{and} \quad \underline{X}_t := \inf_{u \leq t} X_u.$$

The symbol $\mathbf{e}_q$ will always denote a random variable which is independent of $X$ and distributed according to an exponential distribution with parameter $q > 0$. In addition, define, for each $x \in \mathbb{R}$,

$$\tau_x^+ = \inf\{t > 0 : X_t > x\} \quad \text{and} \quad \tau_x^- = \inf\{t > 0 : X_t < x\}.$$

**2. Asymptotic overshoots.** Let us now move to the setting of Klüppelberg, Kyprianou and Maller [10] and, in part, the motivation for this paper. For this it will be necessary to introduce some more notation.

For each $\alpha \geq 0$, $\mathcal{S}^{(\alpha)}$ will denote the class of nonlattice convolution equivalent distributions. That is, to say distributions, $F$, with a nonlattice support on $[0,\infty)$ such that $\overline{F}(x) := 1 - F(x) > 0$ for all $x > 0$ satisfying

(2)
$$\lim_{u \uparrow \infty} \frac{\overline{F}(u-x)}{\overline{F}(u)} = e^{\alpha x} \qquad \text{for each } x \in \mathbb{R},$$

$$\lim_{u \uparrow \infty} \frac{\overline{F}^{*2}(u)}{\overline{F}(u)} = 2M \qquad \text{for some } M > 0.$$

There are several additional facts which follow from this definition. The constant $M$ was identified as equal to $\int_{[0,\infty)} e^{\alpha x} F(dx)$ (and, hence, the latter Laplace–Stieltjes transform is necessarily finite); see [4, 5, 11, 12]. Condition (2) implies that $F(dx)/\overline{F}(x)$ converges in the weak sense to an exponential distribution with parameter $\alpha$. It can also be shown that any measure $\Pi$ which is tail equivalent to a distribution $F \in \mathcal{S}^{(\alpha)}$, that is, to say $\overline{\Pi}(u) := \Pi(u,\infty) \sim \overline{F}(u)$ as $u \uparrow \infty$ for $F \in \mathcal{S}^{(\alpha)}$, also belongs to $\mathcal{S}^{(\alpha)}$; see [7].

The following assumptions are included in the set-up in [10].

ASSUMPTION 1. Fix $\alpha > 0$.

(i) $X_0 = 0$, $\lim_{t \uparrow \infty} X_t = -\infty$ almost surely and $\operatorname{supp} \Pi \cap (0,\infty) \neq \varnothing$,
(ii) $\overline{\Pi}_H \in \mathcal{S}^{(\alpha)}$ and
(iii) $q + \xi(-\alpha) > 0$.

One of the main contributions of Klüppelberg, Kyprianou and Maller [10] was the following result.

THEOREM 2. *Under Assumption 1 we have*
$$\lim_{x \uparrow \infty} P(X_{\tau_x^+} - x > u | \tau_x^+ < \infty) = \overline{G}(u),$$

*where*

(3) $$\overline{G}(u) = \frac{e^{-\alpha u}}{q}\left(q + \xi(-\alpha) + \int_{(u,\infty)} (e^{\alpha y} - e^{\alpha u})\Pi_H(dy)\right).$$

In this paper we aim to recapture and explain in more detail the above result by proving stronger versions of asymptotic results concerning the overshoot and undershoot of both $X$ and $\overline{X}$. Specifically, we will show that the two components

$$\frac{e^{-\alpha u}}{q}(q + \xi(-\alpha)) \quad \text{and} \quad \frac{e^{-\alpha u}}{q}\left(\int_{(u,\infty)} (e^{\alpha y} - e^{\alpha u})\Pi_H(dy)\right)$$

in (3) are the consequence of two types of asymptotic overshoot; namely, first passage occurring as a result of the following:

- an arbitrarily large jump from a finite position after a finite time, or

- a finite jump from a finite distance relative to the barrier after an arbitrarily large time,

respectively.

Our method appeals directly to a new fluctuation identity for a general Lévy process at first passage over a fixed level which specifies the following quintuple law of:

- the time of first passage relative to the time of the last maximum at first passage,
- the time of the last maximum at first passage,
- the overshoot at first passage,
- the undershoot at first passage and
- the undershoot of the the last maximum at first passage.

This quintuple law can be expressed entirely in terms of the quantities $\Pi_X, \mathcal{U}$ and $\widehat{\mathcal{U}}$.

Once this identity is established, it becomes a straightforward exercise to deal with the asymptotic behavior of this quintuple law conditional on first passage occurring under Assumption 1. Indeed, what will prevail in our analysis is the use of the facts that, under this assumption, $U(\cdot, \infty)$ and $\overline{\Pi}_X^+(\cdot)$ both belong to $\mathcal{S}^{(\alpha)}$. These two facts can be deduced from the combined conclusions of Proposition 2.5, Lemma 3.5, Theorem 4.1 and Proposition 5.3 in [10]. Specifically, it was proved that when Assumption 1(i) and (iii) hold, then $U(\cdot, \infty), \overline{\Pi}_H(\cdot)$ and $\overline{\Pi}_X^+(\cdot)$ are all in $\mathcal{S}^{(\alpha)}$ simultaneously or not at all. In the case they all belong to $\mathcal{S}^{(\alpha)}$,

$$(4) \qquad U(u, \infty) \sim \frac{1}{(q + \xi(-\alpha))^2} \overline{\Pi}_H(u) \sim \frac{1}{(q + \xi(-\alpha))^2 \widehat{\xi}(\alpha)} \overline{\Pi}_X^+(u)$$

as $u$ tends to infinity.

The outline of the remainder of the paper is as follows. In the next section we prove the new fluctuation identity for first passage of a general Lévy process over a fixed level. In Section 4 we consider the asymptotic joint laws of the space–time overshoot of $X$, the undershoot of $X$ and the space–time undershoot of $\overline{X}$, all under Assumption 1. We conclude with some additional remarks, in particular, with regard to applications in insurance mathematics.

**3. A quintuple law for overshoots and undershoots.** The main purpose of this section is to prove the following quintuple law for space–time positions of overshoots and undershoots. We will use the notation

$$\overline{G}_t = \sup\{s \leq t : X_s \vee X_{s-} = \overline{X}_{s-}\} \quad \text{and} \quad \underline{G}_t = \sup\{s \leq t : X_s \wedge X_{s-} = \underline{X}_{s-}\}.$$

THEOREM 3. *Suppose that $X$ is not a compound Poisson process. Then for a suitable choice of normalizing constant of the local time at the maximum, for each $x > 0$, we have on $u > 0$, $v \geq y$, $y \in [0, x]$, $s, t \geq 0$,*

$$P(\tau_x^+ - \overline{G}_{\tau_x^+-} \in dt, \overline{G}_{\tau_x^+-} \in ds, X_{\tau_x^+} - x \in du, x - X_{\tau_x^+-} \in dv, x - \overline{X}_{\tau_x^+-} \in dy)$$
$$= \mathcal{U}(x - dy, ds)\widehat{\mathcal{U}}(dv - y, dt)\Pi_X(du + v),$$

where $\Pi_X$ is the Lévy measure of $X$.

Before going to the proof, let us give some intuition behind the statement of this result by discussing its analogue for random walks. The latter turns out to be relatively simple to establish.

Suppose then that $S = \{S_n : n \geq 0\}$ is a random walk on the probability space $(\Omega, \mathcal{F}, \mathbb{P})$. That is, $S_0 = 0$ and $S_n = \sum_{i=1}^n \xi_i$, where $\{\xi_i : i \geq 1\}$ are independent and identically distributed with some law $F$. Define the random variables

$$\overline{S}_n = \max(0, S_1, \ldots, S_n),$$
$$\bar{\theta}^n = \max\{k \leq n : S_k = \overline{S}_n\},$$
$$\sigma_x = \min\{n \geq 1 : S_n > x\}.$$

Let $\{(T'_n, H'_n) : n \geq 0\}$ be the *weak* ascending ladder process and $\{(\widehat{T}_n, \widehat{H}_n) : n \geq 0\}$ be the *strict* descending ladder height process of $S$. Associated with each of these ladder processes are their Green measures

$$\mathbb{U}'(dx, i) := \sum_{n \geq 0} \mathbb{P}(H'_n \in dx, T'_n = i) \quad \text{and} \quad \widehat{\mathbb{U}}(dx, j) = \sum_{n \geq 0} \mathbb{P}(\widehat{H}_n \in dx, \widehat{T}_n = j)$$

for $x \geq 0$ and $i, j \in \mathbb{Z}_{\geq 0}$. The equivalent quintuple law for random walks takes the following form.

THEOREM 4. *For each $x > 0$, we have on $u > 0$, $v \geq y$, $y \in [0, x]$, $i, j \in \mathbb{Z}_{\geq 0}$,*

(5)
$$\mathbb{P}(\sigma_x - 1 - \bar{\theta}^{\sigma_x - 1} = i, \bar{\theta}^{\sigma_x - 1} = j,$$
$$S_{\sigma_x} - x \in du, x - S_{\sigma_x - 1} \in dv, x - \overline{S}_{\sigma_x - 1} \in dy)$$
$$= \mathbb{U}'(x - dy, j)\widehat{\mathbb{U}}(dv - y, i)F(du + v).$$

PROOF. Note first that, by duality,

$$\widehat{\mathbb{U}}(dv - y, i) = \mathbb{P}(S_m < 0, 1 \leq m < i, S_i \in y - dv),$$

so that the

RHS of (5)
$$= \mathbb{P}(S_n < x - y, 1 \leq n < j, S_j \in x - dy)$$
$$\times \mathbb{P}(S_m < 0, 1 \leq m < i, S_i \in y - dv)\mathbb{P}(S_1 \in v + du)$$
$$= \mathbb{P}(S_n < x - y, 1 \leq n < j, S_j \in x - dy,$$
$$S_{j+m} < x - y, 1 \leq m < i, S_{j+i} \in x - dv, S_{j+i+1} \in x + du)$$
$$= \text{LHS of (5)},$$

thus completing the proof. □

REMARK 5. From the analysis above, if we let $\bar{\theta}_n = \min\{k : S_k = \overline{S}_n\}$, then one can reason similarly that, for each $x > 0$, we have on $u > 0$, $v \geq y$, $y \in [0, x]$, $i, j, \in \mathbb{Z}_{\geq 0}$,

$$\mathbb{P}(\sigma_x - 1 - \bar{\theta}_{\sigma_x - 1} = i, \bar{\theta}_{\sigma_x - 1} = j,$$
$$S_{\sigma_x} - x \in du, x - S_{\sigma_x - 1} \in dv, x - \overline{S}_{\sigma_x - 1} \in dy)$$
$$= \mathbb{U}(x - dy, j)\widehat{\mathbb{U}}'(dv - y, i)F(du + v).$$

Here we have the subtle difference that $\mathbb{U}$ and $\widehat{\mathbb{U}}'$ are the Green's measures of the strict ascending and weak descending ladder processes.

Note that hints concerning the quintuple law for the random walk case can already be seen in the discussion on the Wiener–Hopf factorization in [3].

We now move to the proof of the quintuple law for Lévy processes.

PROOF OF THEOREM 3. We prove the result in three steps.

*Step* 1. Let us suppose that $m, k, f, g$ and $h$ are all positive, continuous functions with compact support satisfying $f(0) = g(0) = h(0) = 0$. We prove

(6)
$$E(m(\tau_x^+ - \overline{G}_{\tau_x^+ -})k(\overline{G}_{\tau_x^+ -})f(X_{\tau_x^+} - x)g(x - X_{\tau_x^+ -})h(x - \overline{X}_{\tau_x^+ -}))$$
$$= \widehat{E}_x\biggl(\int_0^{\tau_0^-} m(t - \underline{G}_t)k(\underline{G}_t)h(\underline{X}_t)w(X_t)\,dt\biggr),$$

where $w(z) = g(z)\int_{(z,\infty)} \Pi_X(du)f(u - z)$.

The proof of this result follows by an application of the compensation formula applied to the point process of jumps with intensity measure $dt \times$

$\Pi(dx)$. We have
$$E(m(\tau_x^+ - \overline{G}_{\tau_x^+ -})k(\overline{G}_{\tau_x^+ -})f(X_{\tau_x^+} - x)g(x - X_{\tau_x^+ -})h(x - \overline{X}_{\tau_x^+ -}))$$
$$= E\left(\sum_{t<\infty} m(t - \overline{G}_{t-})k(\overline{G}_{t-})g(x - X_{t-})h(x - \overline{X}_{t-})\right.$$
$$\left. \times \mathbf{1}_{(x-\overline{X}_{t-}>0)} f(X_{t-} + \Delta X_t - x)\mathbf{1}_{(\Delta X_t > x - X_{t-})}\right)$$
$$= E\left(\int_0^\infty dt \cdot m(t - \overline{G}_{t-})k(\overline{G}_{t-})g(x - X_{t-})h(x - \overline{X}_{t-})\right.$$
$$\left. \times \mathbf{1}_{(x-\overline{X}_{t-}>0)} \int_{(x-X_{t-},\infty)} \Pi_X(d\phi) f(X_{t-} + \phi - x)\right)$$
$$= E\left(\int_0^\infty dt \cdot m(t - \overline{G}_{t-})k(\overline{G}_{t-})h(x - \overline{X}_{t-})\mathbf{1}_{(x-\overline{X}_{t-}>0)} w(x - X_{t-})\right)$$
$$= \widehat{E}_x\left(\int_0^\infty dt \cdot \mathbf{1}_{(t<\tau_0^-)} m(t - \underline{G}_t)k(\underline{G}_t)h(\underline{X}_t)w(X_t)\right),$$

which is equal to the right-hand side of (6). In the last equality we have rewritten the previous equality in terms of the path of $-X$. Note that the condition $f(0) = g(0) = h(0) = 0$ has been used implicitly to exclude from the calculation the case of first passage by creeping.

*Step* 2. Next we prove that

$$E_x\left(\int_0^{\tau_0^-} m(t - \underline{G}_t)k(\underline{G}_t)h(\underline{X}_t)w(X_t)\,dt\right)$$
(7)
$$= \int_{[0,\infty)}\int_{[0,\infty)} \mathcal{U}(d\phi, dt)$$
$$\times \int_{[0,x]}\int_{[0,\infty)} \widehat{\mathcal{U}}(d\theta, ds)m(t)k(s)h(x-\theta)w(x+\phi-\theta).$$

(Note, however, that this result will be applied in conjunction with the conclusion of step 1 to the process $-X$.)

The statement and proof of (7) is a generalization of Theorem VI.20 in [2]. For $q > 0$,

$$E_x\left(\int_0^{\tau_0^-} dt \cdot m(t - \underline{G}_t)k(\underline{G}_t)h(\underline{X}_t)w(X_t)e^{-qt}\right)$$
$$= q^{-1}E_x(m(\mathbf{e}_q - \underline{G}_{\mathbf{e}_q})k(\underline{G}_{\mathbf{e}_q})h(\underline{X}_{\mathbf{e}_q})w(X_{\mathbf{e}_q} - \underline{X}_{\mathbf{e}_q} + \underline{X}_{\mathbf{e}_q}); \mathbf{e}_q < \tau_0^-)$$
$$= q^{-1}\int_{[0,x]}\int_{[0,\infty)} P(-\underline{X}_{\mathbf{e}_q} \in d\theta, \underline{G}_{\mathbf{e}_q} \in ds)k(s)$$

$$\times \int_{[0,\infty)} \int_{[0,\infty)} P(X_{\mathbf{e}_q} - \underline{X}_{\mathbf{e}_q} \in d\phi, \mathbf{e}_q - \underline{G}_{\mathbf{e}_q} \in dt)$$

(8)
$$\times m(t)h(x-\theta)w(x+\phi-\theta)$$

$$= q^{-1} \int_{[0,x]} \int_{[0,\infty)} P(-\underline{X}_{\mathbf{e}_q} \in d\theta, \underline{G}_{\mathbf{e}_q} \in ds)k(s)$$

$$\times \int_{[0,\infty)} \int_{[0,\infty)} P(\overline{X}_{\mathbf{e}_q} \in d\phi, \overline{G}_{\mathbf{e}_q} \in dt)$$

$$\times m(t)h(x-\theta)w(x+\phi-\theta),$$

where the Wiener–Hopf factorization and duality have been used in the second and third equalities respectively. Next note that, for a suitable normalization of the local time at the maximum, we have

$$q = \kappa(q,0)\widehat{\kappa}(q,0)$$

(cf. equation (3) of Chapter VI in [2]). Further, it is also known from the Wiener–Hopf factorization that

$$\frac{1}{\kappa(q,0)} E(e^{-\alpha \overline{G}_{\mathbf{e}_q} - \beta \overline{X}_{\mathbf{e}_q}}) = \frac{1}{\kappa(\alpha+q,\beta)}$$

(cf. equation (1), Chapter VI of [2]) and, hence, recalling (1), it follows that

$$\lim_{q \downarrow 0} \frac{1}{\kappa(q,0)} P(\overline{X}_{\mathbf{e}_q} \in d\phi, \overline{G}_{\mathbf{e}_q} \in dt) = \mathcal{U}(d\phi, dt)$$

in the sense of vague convergence. A similar convergence holds for $P(-\underline{X}_{\mathbf{e}_q} \in d\theta, \underline{G}_{\mathbf{e}_q} \in ds)/\widehat{\kappa}(q,0)$. Equality (7) thus follows by taking limits in (8).

*Step* 3. We combine the conclusions of steps 1 and 2 to conclude that

$$E(m(\tau_x^+ - \overline{G}_{\tau_x^+-})k(\overline{G}_{\tau_x^+-})f(X_{\tau_x^+} - x)g(x - X_{\tau_x^+-})h(x - \overline{X}_{\tau_x^+-}))$$

$$= \int_{u>0, y\in[0,x], 0<y\leq v, s\geq 0, t\geq 0} m(t)k(s)f(u)g(v)h(y)$$

$$P(\tau_x^+ - \overline{G}_{\tau_x^+-} \in dt, \overline{G}_{\tau_x^+-} \in ds,$$

$$X_{\tau_x^+} - x \in du, x - X_{\tau_x^+-} \in dv,$$

$$x - \overline{X}_{\tau_x^+-} \in dy)$$

$$= \int_{[0,\infty)} \int_{[0,\infty)} \widehat{\mathcal{U}}(d\phi, dt) \int_{[0,x]} \int_{[0,\infty)} \mathcal{U}(d\theta, ds) m(t)k(s)$$

$$\times h(x-\theta)g(x+\phi-\theta)$$

$$\times \int_{(x+\phi-\theta,\infty)} \Pi_X(d\eta) f(\eta - (x+\phi-\theta)).$$

Substituting $y = x - \theta$, then $y + \phi = v$ and finally, $\eta = v + u$ in the right-hand side above yields

$$E(m(\tau_x^+ - \overline{G}_{\tau_x^+-})k(\overline{G}_{\tau_x^+-})f(X_{\tau_x^+} - x)g(x - X_{\tau_x^+-})h(x - \overline{X}_{\tau_x^+-}))$$
$$= \int_{[0,\infty)}\int_{[0,x]} \mathcal{U}(x - dy, ds) \int_{[0,\infty)}\int_{[y,\infty)} \widehat{\mathcal{U}}(dv - y, dt)$$
$$\times \int_{(0,\infty)} \Pi_X(du + v) m(t) k(s) f(u) g(v) h(y)$$

and the statement of the theorem follows. □

The missing case of a compound Poisson process, excluded from Theorem 3, can be handled similarly to the random walk case.

As a consequence of the above identity, we also obtain the following corollary which relates $\mathbf{\Pi}(dt, dh)$, the Lévy measure of $(L^{-1}, H)$, to $\Pi_X$.

COROLLARY 6. *For all $t, h > 0$, we have*

$$\mathbf{\Pi}(dt, dh) = \int_{[0,\infty)} \widehat{\mathcal{U}}(d\theta, dt) \Pi_X(dh + \theta).$$

PROOF. The result will follow by first proving the following auxiliary result for the ascending ladder process at its first passage time $T_x := \inf\{t > 0 : H_t > x\}$. Let $\Delta L_{T_x}^{-1} = L_{T_x}^{-1} - L_{T_x-}^{-1}$, then

(9)
$$P(\Delta L_{T_x}^{-1} \in dt, L_{T_x-}^{-1} \in ds, x - H_{T_x-} \in dy, H_{T_x} - x \in du)$$
$$= \mathcal{U}(x - dy, ds)\mathbf{\Pi}(dt, du + y)$$

for $t > 0$, $s > 0$, $y \in [0, x]$, $u > 0$. The proof follows from a straightforward calculation using the compensation formula along the lines of the proof of Proposition III.2 in [2]. We omit the technicalities for the sake of brevity.

To finish the proof of the corollary, note that $\Delta L_{T_x}^{-1} = \tau_x^+ - \overline{G}_{\tau_x^+-}$, $L_{T_x-}^{-1} = \overline{G}_{\tau_x^+-}$, $x - H_{T_x-} = x - X_{\tau_x^+-}$ and $H_{T_x} - x = X_{\tau_x^+} - x$. Hence, from the quintuple law, we also know that

$$P(\Delta L_{T_x}^{-1} \in dt, L_{T_x-}^{-1} \in ds, x - H_{T_x-} \in dy, H_{T_x} - x \in du)$$
$$= \mathcal{U}(x - dy, ds) \int_{[y,\infty)} \widehat{\mathcal{U}}(dv - y, dt) \Pi_X(du + v)$$

and from here, by comparing with (9), the statement of the theorem follows. □

Note that, by integrating out $dt$ in the conclusion of the above corollary, we recover the recent identity of Vigon [14] for the Lévy measure of the ascending ladder height process.

We conclude this section with examples of Lévy processes for which new, explicit identities can be obtained. Before doing so we make the remark that there are limited examples of Lévy processes for which the exponents $\kappa$ and $\widehat{\kappa}$ are known explicitly in terms of elementary or special functions. Further, of these known examples, there are no known cases for which the inversion in (1) can be performed to give the bivariate measures $\mathcal{U}$ and $\widehat{\mathcal{U}}$ explicitly. Not surprisingly then, our examples do not explore the quintuple law to its full generality.

EXAMPLE 7 (Strictly stable processes). Suppose that $X$ is a strictly stable process with index $\gamma \in (0,2)$. That is to say, a Lévy process satisfying the scaling property $X_t \stackrel{d}{=} t^{1/\gamma} X_1$ for all $t > 0$. The Lévy measure is given (up to a multiplicative constant) by

$$\Pi_X(dx) = \mathbf{1}_{(x>0)} \frac{c^+}{x^{1+\gamma}}\, dx + \mathbf{1}_{(x<0)} \frac{c^-}{|x|^{1+\gamma}}\, dx,$$

where $c^+$ and $c^-$ are two nonnegative real numbers. To avoid trivialities, we take $c^+ > 0$.

For such processes it is known that the ladder process $H$ is a stable subordinator with index $\gamma\rho$, where $\rho = P(X_1 \geq 0)$ and, hence, up to a multiplicative constant $\kappa(0,\beta) = \beta^{\gamma\rho}$ for $\beta \geq 0$. Similarly, up to a multiplicative constant $\widehat{\kappa}(0,\beta) = \beta^{\gamma(1-\rho)}$. For these facts, the reader is again referred to Bertoin [2].

Inverting (1) when $\alpha = 0$, we find that (up to a multiplicative constant)

$$U(dx) = \frac{x^{\gamma\rho-1}}{\Gamma(\gamma\rho)}\, dx,$$

with a similar identity holding for $\widehat{U}(dx)$ except that $\rho$ is replaced by $1-\rho$. Marginalizing the quintuple law to a triple law, we now obtain a new identity for stable processes. Namely,

$$P(X_{\tau_x^+} - x \in du, x - X_{\tau_x^+-} \in dv, x - \overline{X}_{\tau_x^+-} \in dy)$$
$$= \text{const.} \frac{(x-y)^{\gamma\rho-1}(v-y)^{\gamma(1-\rho)-1}}{(v+u)^{1+\gamma}}\, dy\, dv\, du$$

for $y \in [0,x]$, $v \geq y$ and $u > 0$, where the normalizing constant makes the right-hand side a distribution (note that stable processes do not creep and, hence, there is no atom on the event $\{X_{\tau_x^+} = x\}$ to take care of). Performing a triple integral one may show that the constant is thus equal to

$$\frac{\sin \alpha\rho\pi}{\pi} \frac{\Gamma(\alpha+1)}{\Gamma(\alpha\rho)\Gamma(\alpha(1-\rho))}.$$

EXAMPLE 8 (Spectrally positive processes). In this case, the downward ladder height process is a linear drift with gradient 1 killed at rate $\widehat{q} \geq 0$. For this reason it follows that $\widehat{U}(dx) = e^{-\widehat{q}x}\,dx$. This gives the triple law

$$P(X_{\tau_x^+} - x \in du, x - X_{\tau_x^+-} \in dv, x - \overline{X}_{\tau_x^+-} \in dy)$$
$$= e^{-\widehat{q}(v-y)} U(x-dy) \Pi_X(du+v)\,dv$$

for $y \in [0, x]$, $v \geq y$ and $u > 0$.

The Wiener–Hopf factors for spectrally positive Lévy processes are well understood (cf. Chapter VII in [2]). Indeed, it is known that $\widehat{\kappa}(\alpha, \beta) = \Phi(\alpha) + \beta$, where $\Phi$ is the inverse of the Laplace exponent $\psi(\beta) = \log E(e^{-\beta X_1})$ for $\beta \geq 0$. The identification of $U$ via its Laplace transform in (1) thus simplifies to

$$\text{(10)} \qquad \int_{[0,\infty)} e^{-\beta x} U(dx) = \frac{\beta}{\psi(\beta)}.$$

When, in addition, $X$ has bounded variation and drifts to minus infinity, it is possible to give a more explicit identity for the measure $U$ and, hence, for the above expression. In this case $X$ is the difference of a subordinator and a positive drift of rate $c$ such that $E(X_1) < 0$. It is known then that $\widehat{q} = 0$ and $q = |E(X_1)|$ (see, e.g., Section 6 of [10]). By taking Laplace transforms, we see from (10) that

$$U(dx) = \frac{1}{c} \sum_{n \geq 0} \nu^{*n}(dx),$$

where we understand $\nu^{*0}(dx) = \delta_0(dx)$ and

$$\nu(dx) = \frac{1}{c} \Pi_X(x, \infty)\,dx.$$

[Note that the assumption $E(X_1) < 0$ ensures that $c^{-1} \int_{(0,\infty)} \Pi_X(y, \infty)\,dy < 1$.]

Our triple law now takes the form

$$P(X_{\tau_x^+} - x \in du, x - X_{\tau_x^+-} \in dv, x - \overline{X}_{\tau_x^+-} \in dy)$$
$$= \frac{1}{c} \sum_{n \geq 0} \nu^{*n}(x - dy) \Pi_X(du+v)\,dv$$

for $y \in [0, x]$, $v \geq y$ and $u > 0$.

REMARK 9. The latter example is relevant to insurance mathematics. One may compare against similar results in the papers of Gerber and Shiu [8], Dickson and Drekic [6] and Sun and Yang [13], which concern the classical Cramér–Lundberg process (which in our setting is a spectrally positive compound Poisson process drifting to minus infinity).

**4. The asymptotic role of undershoots in overshoots.** In the following two theorems, we consider the asymptotic overshoot and undershoot in space and time at first passage of $X$, conditional on making first passage, as the barrier tends to infinity. The spatial undershoot is measured, in the first case, backward from the barrier and, in the second case, upward from the origin.

THEOREM 10. *Under Assumption* 1:

(i) *For $t \geq 0$, $y \geq 0$, $v \geq y$ and $u > 0$,*

$$\lim_{x \uparrow \infty} P(\tau_x^+ - \overline{G}_{\tau_x^+ -} \in dt, X_{\tau_x^+} - x \in du,$$

$$x - X_{\tau_x^+ -} \in dv, x - \overline{X}_{\tau_x^+ -} \in dy | \tau_x^+ < \infty)$$

$$= \frac{\alpha}{q} e^{\alpha y} \, dy \cdot \widehat{\mathcal{U}}(dv - y, dt) \Pi_X(du + v).$$

(ii) *For $u > 0$, we have*

$$\int_{v \in (0,\infty)} \lim_{x \uparrow \infty} P(X_{\tau_x^+} - x \in du, x - X_{\tau_x^+ -} \in dv | \tau_x^+ < \infty)$$

$$= \frac{\alpha}{q} \int_0^\infty e^{\alpha y} \Pi_H(du + y) \, dy.$$

PROOF. (i) Starting with the main identity given in Theorem 3, marginalizing out $\overline{G}_{\tau_x^+ -}$ and recalling the Pollaczek–Khintchine identity,

$$P(\tau_x^+ < \infty) = qU(x, \infty)$$

(cf. Proposition 2.5 of [10]), we see that the required asymptotic is equal to

$$\lim_{x \uparrow \infty} \frac{U(x - dy)}{qU(x, \infty)} \widehat{\mathcal{U}}(dv - y, dt) \Pi_X(du + v).$$

Note that $U(\cdot, \infty) \in \mathcal{S}^{(\alpha)}$ by Assumption 1 and so by the associated property of weak convergence, the limit follows.

(ii) In the spirit of the calculations in part (i), we may also marginalize the quintuple law to a bivariate law of the overshoot and undershoot and compute, with the help of weak convergence,

$$\lim_{x \uparrow \infty} P(X_{\tau_x^+} - x \in du, x - X_{\tau_x^+ -} \in dv | \tau_x^+ < \infty)$$

$$= \lim_{x \uparrow \infty} \int_{y \in [0,v]} \frac{U(x - dy)}{qU(x, \infty)} \widehat{U}(dv - y) \Pi_X(du + v)$$

$$= \int_0^\infty dy \cdot \frac{\alpha}{q} e^{\alpha y} \mathbf{1}_{(y \leq v)} \widehat{U}(dv - y) \Pi_X(du + v).$$

Finally integrating with respect to $v \in (0, \infty)$, the result follows by applying Fubini's theorem and invoking Vigon's identity (as a special case of Corollary 6) in the form

$$\Pi_H(du + y) = \int_{[y,\infty)} \widehat{U}(dv - y)\Pi_X(du + v),$$

thus concluding the proof. $\square$

THEOREM 11. *Under Assumption 1:*

(i) *For $s, t \geq 0$, $u > 0$, $\theta \geq 0$ and $\phi \leq \theta$,*

$$\lim_{x\uparrow\infty} P(\tau_x^+ - \overline{G}_{\tau_x^+-} \in dt, \overline{G}_{\tau_x^+-} \in ds,$$

$$X_{\tau_x^+} - x \in du, X_{\tau_x^+-} \in d\phi, \overline{X}_{\tau_x^+-} \in d\theta | \tau_x^+ < \infty)$$

$$= \mathcal{U}(d\theta, ds)\widehat{\mathcal{U}}(\theta - d\phi, dt)\frac{\alpha(q + \xi(-\alpha))^2 \widehat{\xi}(\alpha)}{q} e^{-\alpha(u-\phi)}\, du.$$

(ii) *For $u > 0$,*

$$\int_{\phi \in (0,\infty)} \lim_{x\uparrow\infty} P(X_{\tau_x^+} - x \in du, X_{\tau_x^+-} \in d\phi | \tau_x^+ < \infty)$$

$$= \alpha e^{-\alpha u}\frac{(q + \xi(-\alpha))}{q}\, du.$$

PROOF. (i) With a change of variable in the main identity of Theorem 3, we have

$$P(\tau_x^+ - \overline{G}_{\tau_x^+-} \in dt, \overline{G}_{\tau_x^+-} \in ds,$$

(11) $$\qquad X_{\tau_x^+} - x \in du, X_{\tau_x^+-} \in d\phi, \overline{X}_{\tau_x^+-} \in d\theta | \tau_x^+ < \infty)$$

$$= \mathcal{U}(d\theta, ds)\widehat{\mathcal{U}}(\theta - d\phi, dt)\frac{\Pi_X(du + x - \phi)}{qU(x, \infty)},$$

for $\theta \in [0, x)$ and $\phi \in (-\infty, \theta]$. From (4) and the associated weak convergence, we have that

$$\lim_{x\uparrow\infty} \frac{\Pi_X(du + x - \phi)}{qU(x, \infty)} = \frac{(q + \xi(-\alpha))^2 \widehat{\xi}(\alpha)\alpha}{q} e^{-\alpha(u-\phi)}\, du$$

and the result follows.

(ii) From (11) we have that

$$\lim_{x\uparrow\infty} P(X_{\tau_x^+} - x \in du, X_{\tau_x^+-} \in d\phi | \tau_x^+ < \infty)$$

$$= \lim_{x \uparrow \infty} \int_{\theta \in [\phi,x]} U(d\theta)\widehat{U}(\theta - d\phi)\frac{\Pi_X(du + x - \phi)}{qU(x,\infty)}$$

$$= \int_{\theta \in [\phi,\infty)} U(d\theta)\widehat{U}(\theta - d\phi)\frac{(q+\xi(-\alpha))^2\widehat{\xi}(\alpha)\alpha}{q}e^{-\alpha(u-\phi)}\,du$$

$$= \int_{\theta \in [0,\infty)} U(d\theta)e^{\alpha\theta}\mathbf{1}_{(\theta \geq \phi)}\widehat{U}(\theta - d\phi)e^{-\alpha(\theta-\phi)}$$

$$\times \frac{(q+\xi(-\alpha))^2\widehat{\xi}(\alpha)\alpha}{q}e^{-\alpha u}\,du.$$

Finally integrating out with respect to $\phi \in (-\infty,\theta]$, the second part follows by Fubini's theorem and the facts

$$\alpha \int_0^\infty e^{-\alpha x}\widehat{U}(x)\,dx = \frac{1}{\widehat{\xi}(\alpha)}$$

and

$$\alpha \int_0^\infty e^{\alpha x}U(x)\,dx = \frac{1}{q+\xi(-\alpha)};$$

these are special cases of (1) and the corresponding result for $\widehat{X}$. $\square$

We conclude with some additional remarks following from the results above.

*Asymptotic independence.* Note that in the last theorem we see an intuitively obvious independence appearing between the overshoot and the undershoot.

*Decomposing the law of the asymptotic overshoot.* The conclusions of Theorems 10 and 11 both reprove and provide an interesting explanation for the identity in Theorem 2. A straightforward calculation on the identity in Theorem 10(ii) shows that

$$\int_{v \in (0,\infty)} \lim_{x \uparrow \infty} P(X_{\tau_x^+} - x > u, x - X_{\tau_x^+-} \in dv | \tau_x^+ < \infty)$$

$$= \frac{e^{-\alpha u}}{q}\left\{\int_u^\infty (e^{\alpha y} - e^{\alpha u})\Pi_H(dy)\right\}.$$

Similarly, from Theorem 11(ii), we have

$$\int_{\phi \in (0,\infty)} \lim_{x \uparrow \infty} P(X_{\tau_x^+} - x > u, X_{\tau_x^+-} \in d\phi | \tau_x^+ < \infty) = \frac{e^{-\alpha u}}{q}(q+\xi(-\alpha)).$$

Adding these two identities together recovers the conclusion of Theorem 2. It also shows that the distribution of the conditional asymptotic overshoot has a contribution coming from an arbitrarily large jump at a finite position and after a finite time, or a finite jump from a finite distance relative to the barrier after an arbitrarily large time. Note also from part (i) of the two theorems in this section that, when, asymptotically, the undershoot is close to the barrier, the time of occurrence of the last maximum prior to first passage was historically close to the first passage time. Further, when there is asymptotic first passage due to an arbitrarily large jump, this jump happens early on in the path of the Lévy process.

For further results concerning asymptotic overshoots of Lévy processes (spectrally positive compound Poisson processes) with subexponential tails, see [1].

*Other identities.* There are a number of other identities one can extract from Theorems 10 and 11. For example, one can obtain an expression for the joint law of the asymptotic overshoot of $X$ and undershoot of $H$ measured from the barrier or measured from zero. In the latter case, with an appropriate marginalization of the quintuple law, one easily recovers the identity given in Theorem 4.2(iii) of [10]. This identity says that

$$\lim_{x\uparrow\infty} P(\overline{X}_{\tau_x^+-} \leq z | \tau_x^+ < \infty) = \frac{(q+\xi(-\alpha))^2}{q} \int_{[0,z]} e^{\alpha\theta} U(d\theta).$$

The proof is straightforward and left as an exercise.

*Asymptotic creeping.* From the distribution $G$ given in Theorem 2, one sees that there is an atom at zero of mass $\alpha c/q$. This atom corresponds to the asymptotic conditional probability of creeping over the barrier as it tends to infinity. This can also be derived directly by noting from [9] that, when the drift $c$ of $H$ is positive, $U$ is absolutely continuous and

$$P(X_{\tau_x^+} = x) = cu(x),$$

where $u(x) = dU(x)/dx$. Weak convergence of $U(dx)/U(x,\infty)$ under Assumption 1 now ensures that

$$\lim_{x\uparrow\infty} P(X_{\tau_x^+} = x | \tau_x^+ < \infty) = \lim_{x\uparrow\infty} \frac{cu(x)}{qU(x,\infty)} = \frac{c\alpha}{q}.$$

*Applications to insurance mathematics.* The motivation for the work in [10] came from insurance mathematics and, in particular, the classical ruin problem. The refinements of their results given here also offer direct insight into ruinous behavior.

Within the current context, one may think of $-X$ as the capital of an insurance firm, the so-called risk process. In which case the event of ruin with an initial capital of $x$ units corresponds to the process $X$ starting at the origin and making first passage at $x$. Understanding the conditional asymptotics as $x$ tends to infinity thus gives information about how ruin occurs when the initial revenue of the insurance firm is extremely large.

The classical risk process is the Cramér–Lundberg model which corresponds to $X$ being a spectrally positive compound Poisson process with negative drift. A more suitable generalization, however, corresponds to the case that $X$ is a general spectrally positive Lévy process. In this case, recalling the Lévy–Itô decomposition, one sees a more realistic feature as follows. Large jumps (of magnitude greater than one) correspond to large claims offset by premiums collected at a constant rate corresponding to linear drift. Large jumps occur spaced out by independent exponentially distributed periods of time and, thus, reasonably correspond to disasters. The compensated small jumps which occur with countable, but none the less, unbounded frequency correspond to minor claims; their compensation can be understood as the aggregate of premiums called in to offset the high intensity of claims.

The case that $X$ is spectrally positive also has the advantage that many of the identities given above simplify further. Write $\psi(\theta) = \log E(e^{-\theta X_1})$ for the Laplace exponent. Since the descending ladder height process is nothing more than linear drift, we also have $\widehat{U}(dx) = dx$, $\widehat{\xi}(\alpha) = \alpha$ and $q + \xi(-\alpha) = -\psi(-\alpha)/\alpha$. From the latter, it is also straightforward to deduce that $q = |E(X_1)| < \infty$; see [10] for further details. Our earlier results now tell us, for example, that

$$\lim_{x\uparrow\infty} P(X_{\tau_x^+} - x \in du, x - X_{\tau_x^+ -} \in dv, x - \overline{X}_{\tau_x^+ -} \in dy | \tau_x^+ < \infty)$$
$$= \frac{\alpha}{|E(X_1)|} e^{\alpha y} \, dy \cdot dv \cdot \Pi_X(du + v)$$

for $y \geq 0$, $v \geq y$ and $u > 0$, and

$$\lim_{x\uparrow\infty} P(X_{\tau_x^+} - x \in du, X_{\tau_x^+ -} \in d\phi, \overline{X}_{\tau_x^+ -} \in d\theta | \tau_x^+ < \infty)$$
$$= \frac{\psi(-\alpha)^2}{|E(X_1)|} e^{-\alpha(u-\phi)} U(d\theta) \cdot d\phi \cdot du$$

for $\theta \geq 0$, $\phi \leq \theta$, $u > 0$. Note also that the renewal measure $U$ can now be identified directly in terms of $\psi$, namely,

$$\int_0^\infty e^{-\beta x} U(dx) = \frac{\beta}{\psi(\beta)} \qquad \text{for } \beta > 0.$$

**Acknowledgments.** We should like to thank two anonymous referees whose remarks led to an improved version of this paper. Thanks are also due to Maaike Verloop who showed as part of her "onderzoekopdracht" at Utrecht University that the calculations for Lévy processes may also be imitated for random walks to give the results in Theorem 4 and Remark 5, and to Erik Baurdoux who computed the constant in Example 7. This work was initiated when both authors were visiting ETH Zürich. We are grateful for their hospitality.

## REFERENCES

[1] ASMUSSEN, S. and KLÜPPELBERG, C. (1996). Large deviations results for subexponential tails with applications to insurance risk. *Stoch. Process. Appl.* **64** 103–125. MR1419495

[2] BERTOIN, J. (1996). *Lévy Processes*. Cambridge Univ. Press. MR1406564

[3] BOROVKOV, A. A. (1976). *Stochastic Processes in Queueing Theory*. Springer, New York. MR391297

[4] CHOVER, J., NEY, P. and WAINGER, S. (1973). Degeneracy properties of subcritical branching processes. *Ann. Probab.* **1** 663–673. MR348852

[5] CLINE, D. B. H. (1987). Convolutions of distributions with exponential and subexponential tails. *J. Austral. Math. Soc. Ser. A* **43** 347–365. MR904394

[6] DICKSON, D. C. M. and DREKIC, S. (2004). The joint distribution of the surplus prior to ruin and the deficit at ruin in some Sparre Andersen models. *Insurance Math. Econ.* **34** 97–107. MR2035340

[7] EMBRECHTS, P. and GOLDIE, C. M. (1982). On convolution tails. *Stoch. Process. Appl.* **13** 263–278. MR671036

[8] GERBER, H. U. and SHIU, E. S. W. (1997). The joint distribution of the time of ruin, the surplus immediately before ruin, and the deficit at ruin. *Insurance Math. Econ.* **21** 129–137. MR1604928

[9] KESTEN, H. (1969). *Hitting Probabilities of Single Points for Processes With Stationary Independent Increments*. Amer. Math. Soc., Providence, RI. MR272059

[10] KLÜPPELBERG, C., KYPRIANOU, A. E. and MALLER, R. A. (2004). Ruin probabilities and overshoots for general Lévy insurance risk processes. *Ann. Appl. Probab.* **14** 1766–1801. MR2099651

[11] ROGOZIN, B. A. (2000). On the constant in the definition of subexponential distributions. *Theory Probab. Appl.* **44** 409–412. MR1751486

[12] SHIMURA, T. and WATANABE, T. (2004). Infinite divisibility and generalized subexponentiality. Preprint.

[13] SUN, L. and YANG, H. (2004). On the joint distributions of surplus immediately before ruin and the deficit at ruin for Erlang(2) risk processes. *Insurance Math. Econ.* **34** 121–125. MR2035342

[14] VIGON, V. (2002). Votre Lévy rampe-t-il? *J. London Math. Soc. (2)* **65** 243–256. MR1875147

SCHOOL OF MATHEMATICS
UNIVERSITY OF MANCHESTER
OXFORD ROAD
MANCHESTER M13 9PL
UK
E-MAIL: rad@ma.man.ac.uk
URL: www.maths.man.ac.uk/DeptWeb/Homepages/rad/

SCHOOL OF MATHEMATICS
AND COMPUTER SCIENCE
HERIOT WATT UNIVERSITY
EDINBURGH EH14 4AS
UK
E-MAIL: kyprianou@ma.hw.ac.uk
URL: www.ma.hw.ac.uk/~kyprianou/